\documentclass[a4paper, 12pt]{article}
\usepackage{enumerate, theorem}
\usepackage{amsmath, amsfonts, amssymb}
\usepackage[height=22.5cm, width=15cm]{geometry}
\usepackage[all]{xy}
\usepackage{mathrsfs, yfonts}

\DeclareMathOperator{\id}{id}

 \def\C{{\mathbb C}}

\newcommand{\parag}[1]{\paragraph{\sc{#1.}} }

\newtheorem{thm}{Theorem}[subsection]
\newtheorem{defn}[thm]{Definition}
\newtheorem{cor}[thm]{Corollaire}
\newtheorem{prop}[thm]{Proposition}
\newtheorem{lemma}[thm]{Lemma}

\begin{document}

\date{12/05/15}

\author{Daniel Barlet\footnote{Institut Elie Cartan, Alg\`{e}bre et G\'eom\`{e}trie,\newline
Universit\'e de Lorraine, CNRS UMR 7502   and  Institut Universitaire de France.}.}

\title{Meromorphic quotients for some holomorphic G-actions.}

\maketitle

\parag{Abstract} Using mainly  tools from  [B.13]  and [B.15] we give a necessary and sufficient condition in order that a holomorphic  action of a connected complex Lie group $G$ on a reduced  complex space $X$  admits a strongly quasi-proper meromorphic quotient. We apply this characterization to obtain a result which assert  that, when $G = K.B$ \  with $B$   a closed complex subgroup of $G$ and $K$ a real compact subgroup of $G$, the existence of a strongly quasi-proper meromorphic quotient for the $B-$action  implies, assuming moreover that there exists a $G-$invariant  Zariski open dense subset in $X$ which is good for the $B-$action, the existence of a strongly quasi-proper meromorphic quotient for the $G-$action on $X$.

\parag{AMS classification} 32 M 05, 32 H 04, 32 H 99, 57 S 20.

\parag{Key words} Strongly quasi-proper map, strongly quasi-proper meromorphic quotient, Holomorphic G-action, finite type cycles,

\tableofcontents

\section{Introduction}
In this article we explain how the tools developed in [M.00], [B.08], [B.13]  and [B.15] can be applied to produce in suitable cases a meromorphic quotient of a holomorphic action of a connected complex Lie group $G$ on a reduced  complex space $X$. This uses the notion of {\bf strongly quasi-proper map} introduced in {\it loc. cit.} and our first goal is to give three hypotheses, called ${\rm [H.1], [H.2], [H.3]}$ on the group action which are {\bf equivalent} to the existence of a {\bf strongly quasi-proper meromorphic quotient}, notion defined in the section 1.2.

The proof of this equivalence is the content of proposition \ref{necessary} and theorem \ref{mero. quot.}. Then we discuss these hypotheses and give a simple sufficient condition ${\rm [H.1str]}$ which implies ${\rm [H.1]}$. The existence theorem for a strongly quasi-proper meromorphic quotient under our three assumptions is applied to prove the following results :

\begin{thm}\label{G = K.B}
Assume that we have a holomorphic action of a  connected complex Lie group $G$ on a reduced complex space $X$. Assume that  $G = K.B$ where $K$ is a compact (real) subgroup of $G$ and $B$ a connected complex closed subgroup of $G$. Assume that the action of $B$ on $X$  satisfies   the condition ${\rm [H.1str]}$ on a $G-$invariant Zariski open dense subset $\Omega$ in $X$, and  the conditions ${\rm [H.2]}$ and ${\rm [H.3]}$. Then the $G-$action satisfies ${\rm [H.1str],[H.2]}$ and ${\rm [H.3]}$ ; so it has a strongly quasi-proper meromorphic quotient.
\end{thm}

\parag{Acknowledgements}  An important part of this article comes from discussions with Peter Heinzner during a stay in Bochum. I want to thank him for his help  and for his hospitality.

\section{Strongly quasi-proper meromorphic quotients.}

\subsection{Preliminaries.}

For the definition of the topology on the space $\mathcal{C}_{n}^{f}(X)$ of finite type $n-$cycles in $X$ and its relationship with the (topological) space $\mathcal{C}_{n}^{loc}(X)$ we refer to [B-M],  [B.13] and [B.15].\\
For the convenience of the reader we recall shortly here the definition of a geometrically f-flat map (GF map) and of a strongly quasi-proper map (SQP map) between irreducible complex spaces and we give  a short summary on some properties of the SQP maps. For more details on these notions see [B.13] and [B.15].

\begin{defn}\label{GF} A holomorphic map $f : M \to N$ between two  irreducible irreducible complex spaces is called a {\bf geometrically f-flat map (a GF-map for short)} if the following conditions are fullfilled :
\begin{enumerate}[i)]
\item The map is quasi-proper equidimensionnal and surjective.\\
 Let $n : \dim M - \dim N$.
\item There exists a holomorphic map $\varphi : N \to \mathcal{C}^{f}_{n}(M)$\footnote{that is to say a f-analytic family of finite type $n-$cycles in $M$ parametrized by $N$.} such that for $y$ generic in $N$ the cycle $\varphi(y)$ is reduced and equal to the set-theoretic fiber $f^{-1}(y)$ of $f$ at $y$.
\end{enumerate}
A holomorphic map $f : M \to N$ between two  irreducible irreducible complex spaces will be {\bf strongly quasi-proper (SQP map for short)} if there exists a modification\footnote{a modification is, by definition,  always proper.} 
  $\tau : \tilde{N} \to N$ such that the strict transform\footnote{By definition $\tilde{M}$ is the irreducible component of $M\times_{N}\tilde{N}$ which surjects on $\tilde{N}$ and $\tilde{f}$ is induced by the projection.} $\tilde{f} : \tilde{M} \to \tilde{N}$ of $f$ by $\tau$ is a GF map.\\
  A meromorphic map $M \dasharrow N$ will be called strongly quasi-proper when the projection on $N$ of its  graph is a SQP map.
\end{defn}

Note that a GF map has, by definition, a {\it holomorphic fiber map} and that a SQP  holomorphic (or meromorphic) map has a {\it meromorphic fiber map} via the composition of the holomorphic fiber map of $\tilde{f}$ with the (holomorphic) direct image map $\tau_{*} :  \mathcal{C}^{f}_{n}(\tilde{M}) \to  \mathcal{C}^{f}_{n}(M)$. Of course, a SQP holomorphic map is quasi-proper, but the converse is not true. The notion of strongly quasi-proper map is stable by modification of the target space, property which is not true in general for a quasi-proper map having ``big fibers''  (see [B.15]).\\

Let $\pi : M \to N$ be a SQP map between irreducible complex spaces and define $n := \dim M - \dim N$. By definition of a SQP map, we can find a Zariski open dense subset $N_{0}$ in $N$ and a holomorphic map $\varphi_{0} : N_{0} \to \mathcal{C}_{n}^{f}(M)$ such that 
\begin{enumerate}[i)]
\item For each $y$ in $N_{0}$ we have the equality of subsets $\vert \varphi_{0}(y)\vert = \pi^{-1}(y)$.
\item For $y$ generic in $N_{0}$ the cycle $\varphi_{0}(y)$ is reduced.
\end{enumerate}
Let $\Gamma \subset N_{0}\times \mathcal{C}_{n}^{f}(M)$ be the graph of $\varphi_{0}$. Then also by definition of a SQP map, the closure $\bar \Gamma$ of $\Gamma$ in $N \times \mathcal{C}_{n}^{f}(M)$ is proper over $N$. Then, using the semi-proper direct image theorem 2.3.2 of [B.15], this implies that $\tilde{N} := \bar \Gamma$ is an irreducible complex space (locally of finite dimension) with the structure sheaf induced by the sheaf of holomorphic functions on $N \times \mathcal{C}_{n}^{f}(M)$. Moreover the natural projection $\tau : \tilde{N} \to N$ is a (proper) modification.\\
Let $\tilde{M} := M\times_{N,str} \tilde{N}$ the strict transform of $M$ by $\tau$, that is to say the irreducible component of $M\times_{N}\tilde{N}$ containing the graph of  $\pi_{0}$, the restriction of $\pi$ to the open set $\pi^{-1}(N_{0})$\footnote{this graph is a Zariski open set in $M\times_{N}\tilde{N}$ which is irreducible as $N$ is irreducible.}. Then let $\tilde{\pi} : \tilde{M} \to \tilde{N}$ the strict transform of $\pi$ by the modification  $\tau$, which is induced on $\tilde{M}$  by the natural projection of $M\times_{N}\tilde{N}$ onto $\tilde{N}$. The set-theoretical fiber at $\tilde{y} := (y,C) \in \tilde{N}$ of $\tilde{\pi}$ is the subset $\vert C\vert \times\tilde{y}$ in $\tilde{M}$. The map $\psi : \tilde{N} \to  \mathcal{C}_{n}^{f}(\tilde{M})$ \ given by $(y, C) \mapsto C \times\{\tilde{y}\}$ is holomorphic and satisfies $\vert \psi(\tilde{y})\vert =\tilde{\pi}^{-1}(\tilde{y})$ for all $\tilde{y}$ in $\tilde{N}$. Moreover $\psi(\tilde{y})$ is a reduced cycle for generic $\tilde{y}$ in $ \tilde{N}$. So the map $\tilde{\pi}$ is geometrically f-flat. It is the {\bf canonical GF-flattning} of $\pi$.\\
Then we have an isomorphism induced by $\psi$
$$ \psi : \tilde{N} \longrightarrow \mathcal{C}_{n}^{f}(\tilde{\pi})  $$
where $ \mathcal{C}_{n}^{f}(\tilde{\pi})      := \{ C \in  \mathcal{C}_{n}^{f}(\tilde{M}) \ / \  \exists \tilde{y} \in \tilde{N} \ s. t. \ \vert C\vert \subset \tilde{\pi}^{-1}(\tilde{y})\}$ is a closed analytic subset of $ \mathcal{C}_{n}^{f}(\tilde{M})$ (see [B.15] proposition 2.1.7.);  the inverse map is induced by the holomorphic map $\hat{\pi} : \mathcal{C}_{n}^{f}(\tilde{\pi}) \to \tilde{N}$ which associates to $\gamma \in \mathcal{C}_{n}^{f}(\tilde{\pi}) $ the point in $\tilde{N}$ whose $\tilde{\pi}-$fiber contains $\gamma$ (see the proposition 2.1.7 of [B.15]).\\
The direct image of $n-$cycles by $\tau$ gives a holomorphic map $\tau_{*} :  \mathcal{C}_{n}^{f}(\tilde{M}) \to \mathcal{C}_{n}^{f}(M)$ which sends $\tilde{N} \simeq \mathcal{C}_{n}^{f}(\tilde{\pi})$ in $\mathcal{C}_{n}^{f}(\pi)$. Let us show that it is an isomorphism of $\tilde{N}$ onto its image in $\mathcal{C}_{n}^{f}(\pi)$ :\\
We have an obvious holomorphic map $\tilde{N} \to \mathcal{C}_{n}^{f}(\pi)$ given by $(y,C) \mapsto C$. We have also a holomorphic map $\mathcal{C}_{n}^{f}(\pi) \to N \times \mathcal{C}_{n}^{f}(M)$ given by $C \mapsto (\hat{\pi}, C)$ where $\hat{\pi} : \mathcal{C}_{n}^{f}(\pi) \to N$ is the map associating to $C \in \mathcal{C}_{n}^{f}(\pi) $ the point $y \in N$ such that $\vert C\vert \subset \pi^{-1}(y)$. This proves our claim.\\
Remark that $\mathcal{C}_{n}^{f}(\pi) $ is not, in general, a complex space (locally of finite dimension).\\

\subsection{Action of $G$ on $\mathcal{C}_{n}^{f}(X)$.}

Let $G$ be a Lie group. We shall say that $G$ acts {\bf continuously holomorphically} on the reduced complex space $X$ when the action $f : G \times X \to X$ is a continuous map such that for each $g \in G$ fixed, the map $x \mapsto f(g,x)$ is a (biholomorphic)  automorphism of $X$. Then there is a natural action of $G$ induced on the set $\mathcal{C}_{n}^{f}(X)$ of finite type $n-$cycles given by $(g,C) \mapsto g_{*}(C)$ where we denote $g_{*}(C)$ the direct image of the cycle $C$ by the automorphism of $X$ associated to $g \in G$. When $G$ is a complex Lie group and the map $f$ is holomorphic  we shall say that the action is {\bf completely holomorphic}.
 
 \begin{prop}\label{holomorphy}
 The action of $G$ on $\mathcal{C}_{n}^{f}(X)$  is continuously holomorphic. This precisely means that the map $G \times\mathcal{C}_{n}^{f}(X) \to \mathcal{C}_{n}^{f}(X)$ given by  $(g,C) \mapsto g_{*}(C)$  is continuous and that, for any f-analytic family of $n-$cycles $(C_{s})_{s \in S}$ in $X$ parametrized by a reduced complex space $S$, the family
 $g_{*}(C_{s})_{(g, s) \in \{g\}\times S}$ is f-analytic for each fixed $g\in G$ . If $G$ is complex Lie group and the action is completely holomorphic, the action of  $G$ on $\mathcal{C}_{n}^{f}(X)$ is completely holomorphic, so for  any f-analytic family of $n-$cycles $(C_{s})_{s \in S}$ in $X$ parametrized by a reduced complex space $S$, the family
 $g_{*}(C_{s})_{(g, s) \in G\times S}$ is f-analytic.
 \end{prop}
 
 \parag{Proof} First we prove the continuity of the action of $G$ on $\mathcal{C}_{n}^{loc}(X)$. To apply the theorem IV 2.5.6 de [B-M] it is enough to see that the map $ F : G \times X \to G\times X$ given by $(g, x) \mapsto (g, g.x)$ is proper. But if $L \subset G$ and $K \subset X$ are compacts sets, we have $F^{-1}(L\times K) \subset L \times (L^{-1}.K)$ which is a compact set in $G \times X$.\\
  The only point left to prove the continuity statement for the topology of $\mathcal{C}_{n}^{f}(X)$, assuming that the continuity for the topology of $\mathcal{C}_{n}^{loc}(X)$ is obtained as follows :\\
 Let $W$ be a relatively compact open set in $X$ and $\mathcal{W}$ be the open set in $\mathcal{C}_{n}^{f}(X)$ of cycles $C$ such any irreducible component of $C$ meets $W$. Then we want to show that the set of $(g, s) \in G \times S$ such that $g_{*} (C_{s})$ lies in $\mathcal{W}$ is an open set in $G \times S$. As the topology of $\mathcal{C}_{n}^{f}(X)$ has a countable basis\footnote{This is a corollary of the fact that this is true for $\mathcal{C}_{n}^{loc}(X)$ (see [B-M] ch.IV) as the topology of $X$ has a countable basis of open sets; see  [B.15] for details.} it is enough to show that if a sequence $(g_{\nu}, s_{\nu})$ converges to $(g, s)$ with $g_{*}(C_{s}) \in \mathcal{W}$ then for $\nu \gg 1$ we have also $(g_{\nu})_{*}(C_{s_{\nu}})\in \mathcal{W}$. If this not the case, we can choose for infinitely many $\nu$ an irreducible component  $\Gamma_{\nu}$ of $(g_{\nu})_{*}(C_{s_{\nu}})$ which does not meet $W$.  Up to pass to a sub-sequence, we may assume that the sequence $\Gamma_{\nu}$ converges in $\mathcal{C}_{n}^{loc}(X)$ to a cycle $\Gamma$ which does not meet $W$ and is contained in $g_{*}(C_{s})$. This is a simple consequence of the continuity of the $G-$action on $\mathcal{C}^{loc}_{n}(X)$ and the characterization of compact subsets in $\mathcal{C}^{loc}_{n}(X)$ (see [B-M] ch.IV). As any irreducible component of $g_{*}(C_{s})$ meets $W$ this implies that $\Gamma$ is the empty $n-$cycle. This means that for any compact $K$ in $X$ there exists an integer $\nu(K)$ such that for $\nu \geq \nu(K)$ we have $\Gamma_{\nu} \cap K = \emptyset$. Choose now a compact neighbourhood $L$ of $g(K)$. For $\nu$ large enough we shall have $K\subset g_{\nu}^{-1}(L)$. This comes from the fact that the automorphisms $g_{\nu}^{-1}$ converge to $g^{-1}$ in the compact-open topology. Then this implies that for $\nu \geq \nu(L)$ the irreducible component $g_{\nu}^{-1}(\Gamma_{\nu})$ of $C_{s_{\nu}}$ does not meet $K$. Then, when $s_{\nu} \to s$ the cycles $C_{s_{\nu}}$ does not converge to $C_{s}$ for the topology of $\mathcal{C}_{n}^{f}(X)$ because we have some ``escape at infinity'' in  a well choosen sub-sequence. Contradiction.$\hfill \blacksquare$\\
 
\begin{lemma}\label{action}
Assume now that in the situation above we have a Lie group $H$ acting continuously holomorphically on $M$. Assume that $\pi$ satisfies $\pi(h.x) = \pi(x)$ for each $x \in M$. Then the strict transform $\tilde{M}$ of $M$ by the canonical modification $\tau : \tilde{N} \to N$ giving the canonical GF flattning of $\pi$ has a natural continuous holomorphic action of $H$ such that the projection $\tilde{M} \to M$ is $H-$equivariant. Moreover the GF map $\tilde{\pi} : \tilde{M} \to \tilde{N}$ satisfies $\tilde{\pi}(h.\tilde{x}) = \tilde{\pi}(\tilde{x})$ for each $(h,\tilde{x}) \in H \times \tilde{M}$.
\end{lemma}

\parag{proof} We have a natural action of $H$ on $\mathcal{C}_{n}^{f}(M)$ which is continuous and holomorphic for each fixed $h \in H$. As the action of $H$ is trivial on the set of fibers of $\pi$, the action of $H$ is trivial on $\tilde{N} $. Define the action of $H$ on $M\times_{N}\tilde{N}$ by the formula
$$ h.(x,(\pi(x), C)) = (h.x, \pi(x), h_{*}(C)) .$$
It is easy to see that this is a continuous holomorphic action, and that it leaves $\tilde{M}$ globally invariant. Now we have for each $h \in H$ and each $(x,\pi(x), C) \in \tilde{M}$ :
$$ \tilde{\pi}\big(h.(x, (\pi(x), C))\big) = \tilde{\pi}\big(h.x, (\pi(x), h_{*}(C))\big) = (\pi(x), h_{*}(C)) = (\pi(x), C) = \tilde{\pi}(x, (\pi(x), C)) $$
because a limit of $H-$invariant cycles is an $H-$invariant cycle.$\hfill \blacksquare$\\

We shall also use the following simple tool from the cycle's space.

\begin{prop}\label{non inf. br.}

Let \ $M$ \ be a reduced complex space and \ $(X_{s})_{s\in S}$ \ a f-continuous family of \ $d-$dimensional finite type cycles parametrized by a  compact subset \ $S$ in $\mathcal{C}_{d}^{f}(M)$. Let \ $(C_{t})_{t \in T}$ \  be a f-continuous family  of finite type non empty  \ $n-$dimensional cycles  in $M$ parametrized by a subset $T$  in $\mathcal{C}_{n}^{f}(M)$ which is compact in \ $\mathcal{C}_{n}^{loc}(M)$. We assume the following condition :
\begin{itemize}
\item  There exists an open dense set $T'$  in $T$ such that  each \ $C_{t}, t \in T'$ \ is equal to the union some \ $X_{s}$. \hfill $(@@)$
\end{itemize}
Then $T$ is a compact subset in $\mathcal{C}^{f}_{n}(M)$.
\end{prop}

\parag{Proof} First remark that, as \ $S$ \ is compact, there exists a compact set \ $L \subset M$ \ such that any irreducible component of any \ $X_{s}$ \ meets \ $L$.\\
Let $(t_{m})_{m \in \mathbb{N}}$ be a  sequence of points in $T'$ converging to a point $t \in T$ and denote by $C_{m}$ the cycle $C_{t_{m}}$ for short and $C_{t} = C_{\infty}$.
Now choose for each \ $m$ \ an irreducible component \ $\Gamma_{m}$ \ of some \ $X_{s_{m}}$ \ contained in \ $C_{m}$. Up to pass to a subsequence, we may assume that \ $\Gamma_{m}$ \ converges in \ $\mathcal{C}_{d}^{f}(M)$ \ to a cycle \ $\Gamma$ \ which is non empty (it contains a point in \ $L$) \ and included in \ $\vert C_{\infty}\vert$. So  \ $C_{\infty}$ \ is not the empty cycle. \\
Let \ $x$ \ be a generic point of an irreducible component \ $D$ \ of \ $C_{\infty}$. Then, up to pass to a subsequence, we may  choose a sequence \ $(x_{m})$ \ of points respectively in \ $C_{m}$ \ which converges to \ $x$. Choose for each \ $m$ \ an irreducible component \ $\Gamma_{m}$ \ of some \ $X_{s_{m}} \subset \vert C_{m}\vert$ \ which contains \ $x_{m}$. This is possible because of condition \ $(@@)$. Now, again up to pass to a subsequence, we may assume that the sequence \ $(\Gamma_{m})_{m \in \mathbb{N}}$ \ converges in \ $\mathcal{C}_{d}^{f}(M)$ \ to a cycle \ $\Gamma$ \ containing the point \ $x$ \ and contained in \ $\vert C_{\infty}\vert$. Note that $\vert \Gamma\vert$ is contained in some $\vert X_{s_{\infty}}\vert$ as we may assume, by compactness of $S$,  that  the sequence $(s_{m})$ converges to $s_{\infty}\in S$. Then we have $\vert X_{s_{\infty}}\vert \subset \vert C_{\infty}\vert$. As \ $D$ \ is the only irreducible component of \ $C_{\infty}$ \ containing \ $x$, it contains at least an irreducible component of  $\vert X_{s_{\infty}}\vert$ containing $x$,  and so \ $D$ \ meets \ $L$. So we have proved that \ $C_{\infty}$ \ is not the the empty \ $n-$cycle and that any irreducible component of \ $C_{\infty}$ \ meets the compact set \ $L$. This is enough to conclude thanks to the proposition 3.2.2 in [B.15]. $\hfill \blacksquare$\\

\subsection{Definition of SQP meromorphic quotient.}

We shall consider a complex connected Lie group $G$ and a completely holomorphic action of $G$ on an irreducible complex space $X$. It is given, by definition, by a holomorphic map $f : G\times X \to X$
such that $f(g.g',x) = f(g, f(g'.x))$ for all $g, g' \in G$ and $x \in X$, assuming that for each $g \in G$ the holomorphic map $x \mapsto f(g,x)$ is an automorphism of $X$, and that $f(1,x) = x$ for all $x \in X$.

\begin{def}\label{SQP-mero-quot.}
A {\bf strongly quasi-proper meromorphic quotient} (we shall say a {\bf SQP-meromorphic quotient} for short) for such an action $f : G\times X \to X$ will be the following data:
\begin{enumerate}
\item a $G-$modification\footnote{This means that we have a completely  holomorphic $G-$action on $\tilde{X}$ and that the modification $\tau$ is $G-$equivariant.} $\tau : \tilde{X} \to X$ with center $\Sigma$.
\item a holomorphic $G-$invariant GF map \ $q : \tilde{X} \to Q$ \ where $Q$ is an irreducible complex space.
\item an analytic $G-$invariant  subset $Y \subset X$ containing $\Sigma$, with no interior point in $X$. We shall denote $\tilde{Y} := \tau^{-1}(Y)$, $\Omega := X \setminus Y, \tilde{\Omega} :=\tau^{-1}(\Omega)$ and $Q' := q(\tilde{\Omega})$. Note that, as $q$ is an open surjective map, $Q'$ is open and dense in $Q$.
\end{enumerate}
Now we ask that these data satisfy the following properties :
\begin{enumerate}[i)]
\item The restriction to $ \Omega$ of the map $q\circ \tau^{-1}$ is a GF map onto the dense open set  $Q'$ in $Q$ and there is an open dense set $Q''$ in $Q'$ such that each fiber  of $q\circ \tau^{-1}$  at a  point in $Q''$ is equal to a $G-$orbit  in $\Omega$.
\item There exists an open dense subset $\Omega_{0} \subset \Omega$ such that for each $\tilde{x}$ in $ \tilde{\Omega}_{0} := \tau^{-1}(\Omega_{0})$ the {\bf closure} $\overline{G.\tilde{x}}$  of  $G.\tilde{x}$ in $\tilde{X}$ is exactly the set $q^{-1}(q(\tilde{x}))$. 
\end{enumerate}
\end{def}

\begin{prop}\label{universal}
Let $G$ be a complex connected Lie group which acts completely holomorphically on an irreducible complex space $X$. Assume that we have a SQP meromorphic quotient for this action, given by a modification $\tau : \tilde{X} \to X$ and a $G-$invariant  GF map $q : \tilde{X} \to Q$. Then let $\psi : Q \to \mathcal{C}_{n}^{f}(X)$ be the holomorphic map obtained by the composition  of the fiber map of the GF map $q$ and the direct image map for $n-$cycles by the modification $\tau$. Define $Q_{u}:= \psi(Q)$. Then we have the following properties:
\begin{enumerate}
\item $Q_{u}$ is a closed analytic subset in $  \mathcal{C}_{n}^{f}(X)$ which is an irreducible complex space (locally of finite dimension) with the structure sheaf induced by the sheaf of holomorphic functions on $ \mathcal{C}_{n}^{f}(X)$.
\item Let $\tilde{X}_{u}$ be the graph of the meromorphic map $ q_{u}: X \dasharrow Q_{u}$ given by the holomorphic  map  $\psi\circ q : \tilde{X} \to Q_{u}$ and let $\tau_{u}: \tilde{X}_{u} \to X$ and $q_{u} : \tilde{X}_{u} \to Q_{u}$  be the projections on $X$ and $Q_{u}$ respectively  of this graph. Then $(\tau_{u}, q_{u})$ is also a SQP meromorphic quotient for the given $G-$action.
\item For any  SQP meromorphic quotient $(\tau, q)$  there exists a unique holomorphic surjective map $\eta : Q \to Q_{u} $ such that the meromorphic maps $q : X \dasharrow Q$ and $q_{u} : X \dasharrow Q_{u}$ 
satisfies $\eta\circ q = q_{u}$.
\item For any G-invariant holomorphic map $h : X \to Y$ there exists a holomorphic map $H : Q_{u} \to Y$ such that  $h\circ\tau_{u} = H \circ q_{u}$.
\end{enumerate}
\end{prop}

\bigskip

\begin{defn}\label{univ.}
In the situation of the previous theorem the SQP meromorphic quotient for the given $G-$action defined by  $(\tau_{u}, q_{u})$  will be called {\bf the minimal SQP meromorphic quotient} of this $G-$action.
\end{defn}

So the proposition above says that the existence of a SQP meromorphic quotient for the given $G-$action implies the existence and uniqueness of a minimal meromorphic quotient for this $G-$action.

\parag{Proof} To prove the point 1) we shall prove that the map $\psi\circ q : \tilde{X} \to  \mathcal{C}_{n}^{f}(X)$ is semi-proper.\\
Let $C \not= \emptyset$ be in $\mathcal{C}_{n}^{f}(X)$ and fix a relatively compact open set $W$ in $X$ meeting all irreducible components of $C$. The subset $\mathcal{W}$ of $\mathcal{C}_{n}^{f}(X)$ of cycles $C'$  such that any irreducible component of $C'$ meets $W$ is an open set containing $C$. Now $ q(\tau^{-1}(\bar W))$ is a compact set in $Q$, as $\tau$ is proper. Take any $y \in Q$ such that $C' : =   \psi(y)$ is in $\mathcal{W}$. The point $y$ is the limit in $Q$ of points $y_{\nu} \in q(\tilde{\Omega}_{0})$ such that the fiber of $Q$ at $y$ is limit in $\mathcal{C}_{n}^{f}(\tilde{X})$ of the fibers $q^{-1}(y_{\nu}) = \overline{G.\tilde{x_{\nu}}}$ where, for $\nu \gg 1$, we can choose $\tilde{x}_{\nu}$ in $\tilde{\Omega}_{0} \cap \tau^{-1}(W)$. Up to pass to a sub-sequence, we may assume that $\tilde{x}_{\nu}$ converges to a point $\tilde{x}$ in  $\tau^{-1}(\bar W)$. Then the continuity of $q$ implies that $q(\tilde{x}) = y$ and $C'$ is the limit of $\overline{G.\tilde{x}_{\nu}}$. So $\vert C'\vert$ is in  the image by $\psi$ of the compact set $q(\tau^{-1}(\bar W))$ and this  gives the semi-properness of $\psi\circ q$.\\
Now the  direct image theorem 2.3.2 in [B.15] shows that $Q_{u}$ is an irreducible complex space (locally of finite dimension) and the point 1) is proved.\\
To prove the second point we have to show that the map $q_{u} : \tilde{X}_{u} \to Q_{u}$ is a GF map. By definition $ \tilde{X}_{u}$ is the closure in $X \times Q_{u}$ of the graph of the map $q_{\vert \Omega_{0}}$ where $\Omega_{0}$ is an open dense set in $X$ such that for any point $x \in \Omega_{0}$ we have $\psi(q(x)) = \overline{G.x}$ (as $\Omega_{0}$ is disjoint from the center of the modification $\tau$ we identify here $\Omega_{0}$ and $\tilde{\Omega}_{0} := \tau^{-1}(\Omega_{0})$). Then by irreducibility of $Q_{u}$ and $X$ the closed analytic subset $ \tilde{X}_{u} \subset X \times Q_{u}$ is equal to the graph of the tautological family of cycles in $X$ parametrized by $Q_{u} \subset \mathcal{C}_{n}^{f}(X)$. This proves the point 2).\\
Now consider a  SQP meromorphic quotient of the given action given by the maps $\tau : \tilde{X} \to X$ and $q : \tilde{X} \to Q$. Let $\psi : \tilde{X} \to \mathcal{C}_{n}^{f}(X)$ the composition of the holomorphic map classifying the fibers of the GF map $q$ with the direct image of $n-$cycles by the modification $\tau$. Then, by the construction in the proof of the point 1), we know that $\psi(Q) = Q_{u}$ and then the map
$\psi$ induces a surjective holomorphic map $\eta : Q \to Q_{u}$. We may  assume that $\tilde{X}$ is in fact the graph of the meromorphic map $q : X \dasharrow Q$. Then, as $\tilde{X}_{u}$ is the graph of the meromorphic map $q_{u} : X \dasharrow Q_{u}$ the holomorphic map $\id_{X}\times \eta : X \times Q \to  X \times Q_{u}$ sends $\tilde{X}$ to $\tilde{X}_{u}$ because this is true over a dense open set in $X$ where the maps $q$ and $q_{u}$ are holomorphic and satisfy $q^{-1}(q(x)) = q_{u}^{-1}(q_{u}(x)) = \overline{G.x}$. This complete the proof of  3).\\
Consider now a G-invariant holomorphic map $h : X \to Y$. Recall (see [B.15] proposition 2.1.7) that $\mathcal{C}_{n}^{f}(h)$, the subset of cycles $C \in \mathcal{C}_{n}^{f}(X)$ which are contained in a fiber of $h$ is a closed analytic subset in $\mathcal{C}_{n}^{f}(X)$. On an open dense subset of $Q_{u}$ we know that a point corresponds to a cycle with support $\overline{G.x}$ for some $x \in X$. Then this means that an open dense subset in $Q_{u}$ is contained in $\mathcal{C}_{n}^{f}(h)$. As this subset is closed we obtain $Q_{u}\subset \mathcal{C}_{n}^{f}(h)$. Now there is a holomorphic map $\hat{h} : \mathcal{C}_{n}^{f}(h) \to Y$ associating to a cycle $C$  the point in $Y$ such $\vert C\vert \subset h^{-1}(y)$ (see [B.15] proposition 2.1.7). This induces a holomorphic map $ H : Q_{u} \to Y$, and it is clear that the relation $h\circ\tau_{u} = H \circ q_{u}$ is true on a dense open set, so everywhere.$\hfill \blacksquare$

\subsection{Good points, good open set.}

Let \ $X$ \ be a reduced complex space and \ $G$ \ be a connected  complex Lie group.
Let \ $f : G \times X \to X$ \ be a holomorphic action of \ $G$ \ on \ $X$. 

 \begin{defn}\label{good}
 We shall say that a point \ $x \in X$ \ \ is a { \bf good point}  for the action \ $f$ \ if the following condition is satisfied
 \begin{itemize}
 \item For each compact set \ $K$ \ in \ $X$ \ there exists an open  neighbourhood \ $V$ \ of \ $x$ \ and a compact set \ $L$ \  in \ $G$ \ such that if \ $y \in V$ \  and \ $g \in G$ \ are such that \ $g.y \in K$,  there exists \ $\gamma \in L$ \ with \ $\gamma.y = g.y$
 \end{itemize}
 We shall say  that the action of $G$ on $X$ is {\bf good} when each point in $X$ is a good point. If $\Omega$ is  a $G-$invariant open set in  $ X$, we shall say that $\Omega$ is a {\bf a good open set }  for the action \ $f$ \  when all  points in $\Omega$ are good points for the $G-$action given by  \ $f$ {\bf restricted to $\Omega$}.\\
 \end{defn}
 
\parag{Remarks}\begin{enumerate}
\item If \ $x \in X$ \ is a good point, then for any \ $g_0 \in G$ \ $g_0.x$ \ is also a good point : for \ $K$ \ given, choose \ $g_0.V$ \ as neighbourhood of \ $g_0.x$ \ and  the compact set \ $L.g_0^{-1} \subset G$ \  to satisfy the needed conditions.
\item If \ $\Omega$ \ is a  good open set, then points in $\Omega$ are not in general good points for the action on $X$.
\item If \ $\Omega$ \ is a  good open set and $W \subset \Omega$ a G-invariant open set, then $W$ is a good open set.
\item If $M$ is a compact set of good points in $X$ for any compact set $K$ in $X$ we can find a neighbourhood $V$ of $M$ in $X$ and a compact set $L$ in $G$ such that for any point $y \in V$ and any $g \in G$ such that $g.y \in K$ there exists $\gamma \in L$ with $\gamma.y = g.y$. This is easily obtained by a standard compactness argument. We shall say that a  compact set of good points is {\bf uniformely good}.\\
\end{enumerate}

\begin{lemma}\label{Good 1}
Let $x$ be a point in $X$. Then $x$ is a good point for the $G-$action on $X$ if and only if the map $F_{X} : G \times X \to X \times X$ given by $(g,x) \mapsto (x, g.x)$ is semi-proper at each point of $\{x\}\times X$. As a consequence  a $G-$invariant open set $\Omega$ in $X$ is a good open set for the $G-$action if and only if the map $ F_{\Omega} : G \times \Omega \to \Omega \times \Omega $ given by $(g, x) \mapsto (x, g.x)$ is semi-proper.
 \end{lemma}
 
 \parag{Proof} Let $x \in X$ be a good point and fix any $z \in X$. To prove that the map $F_{X}$ is semi-proper at  $(x, z)$ choose compact neighbourhoods $V_{0}$ and $K$ of $x$ and $z$ in $X$ and apply the definition of a good point to the compact set $K$. So we can find a neighbourhood $V$ of $x$, that we may assume to be contained in $V_{0}$, and a compact set $L$ in $G$ such that for any $y \in V$ such that $g.y\in K$ we have a $\gamma \in L$ with $g.y = \gamma.y$. Then we have $F_{X}(G \times X) \cap (V\times K) = F_{X}(L \times V_{0}) \cap (V\times K) $ and $L \times V_{0}$ is a compact set in $G \times X$.\\
 Conversely, assume that the map $F_{X}$ is semi-proper at each point of $ \{x\}\times X$. Take a compact set $K$ in $X$ and apply the semi-properness to each point $(x, z)$ where $z$ is in $K$. For each $z \in K$ we obtain  open neighbourhoods $V_{z}$ and $W_{z}$  of $x$ and $z$ in $X$ and a compact set $L_{z}\times M_{z}$ in $G\times X$ such that 
 $$F_{X}(G \times X)\cap (V_{z}\times W_{z}) = F_{X}(L_{z}\times M_{z}) \cap (V_{z}\times W_{z}) .$$
 Extract a finite sub-cover $W_{1}, \dots, W_{N}$  of $K$ by the open sets $W_{z}$ and define the compact set  $L : = \cup_{i \in [1, N]} L_{z_{i}}$ and the neighbourhood  $V := \cap_{i \in [1,N]} V_{z_{i}}$ of $x$. Then if $y$ is in $V$  and $g.y$ is in $K$ there exists $i \in [1,N]$ such that $g.y \in W_{i}$. As $y$ is in $V_{i}$ we can find a $\gamma \in L_{z_{i}}\subset L$ with $F_{X}(g, y) = F_{X}(\gamma, y)$ and this implies that $x$ is a good point. \\
 The second assertion is an easy consequence of the first one.$\hfill \blacksquare$\\

\begin{prop}\label{Good 2}
Let \ $G$ \ be a connected  complex Lie group. Let \ $f : G \times X \to X$ \ be an holomorphic action of \ $G$ \ on a reduced complex space \ $X$. Then we have the following properties :
\begin{enumerate}[i)]
\item If \ $x$ \ is a good point for \ $f$ \ the orbit \ $G.x$ \ is a closed analytic subset of \ $X$.
\item If \ $x$ \ is a good point for \ $f$ \ and if \ $(G.x) \cap K = \emptyset $ \ where \ $K \subset X$ \ is a compact set, there exists a neighbourhood \ $V$ \ of \ $x$ \ in \ $X$ \ such that \ $(\overline{ G.x'}) \cap K = \emptyset$ \ for any \ $x'$ \ in \ $V$ ($x'$ is not assume here to be a good point).
\item If \ $x$ \ is a good point for \ $f$ \ there exists a neighbourhood \ $V$ \ of \ $x$ \ such that any good point \ $x' \in V$ \ has an orbit which is a closed analytic subset of the same dimension than \ $G.x$.
\item Let \ $\Omega$ \ be   a good connected open set  in $X$ which is normal  and let \ $n$ \ be the dimension of \ $G.x$ \ for \ $x \in \Omega$. Then there exists a holomorphic map\footnote{This means that we have an $f-$analytic family of $n-$cycles in \ $\Omega$ \ parametrized by \ $\Omega$.} $\varphi : \Omega \to \mathcal{C}_n^f(\Omega) $   given generically by \ $\varphi(x) : = G.x$ \ as a reduced $n-$cycle in \ $\Omega$.
\item When we have a good open set $\Omega$  in $X$ which is normal, there exists a quasi-proper equidimensional holomorphic quotient of $\Omega$ for the action restricted to $\Omega$. 
\end{enumerate}
\end{prop}

\parag{Proof} We already proved that $x$ is a good point if the map $G \to X$ given by $g \mapsto g.x$ s semi-proper in lemma 1.3.2. Now Kuhlmann's theorem  [K.64], [K.66]  gives that \ $f_x(G) = G.x$ \ is a closed analytic subset of \ $X$. This proved i)\\

Assume ii) is not true ;  then we have a compact set \ $K$ \ such that \ $(G.x)\cap K = \emptyset$ \ and a sequence \ $(x_{\nu})_{\nu \in \mathbb{N}}$ \ converging to \ $x$ \ and such that \ $(\overline{G.x_{\nu}}) \cap K$ \ is not empty for each \ $\nu$. Fix a compact neighbourhood $\tilde{K}$ of $K$ such that $(G.x) \cap \tilde{K} = \emptyset$. This is possible thanks to i). Pick a point \ $y_{\nu} = \lim_{\alpha \to \infty} \ g_{\nu,\alpha}.x_{\nu}$ \ in \ $(\overline{G.x_{\nu}}) \cap K$ \ for each \ $\nu$. Up to pass to a subsequence we may assume that sequence\ $(y_{\nu})$ \ converges to \ $y \in K$ \ when \ $\nu \to + \infty$. So, for $\alpha \geq \alpha(\nu)$, we can assume that $g_{\nu,\alpha}.x_{\nu}$ is in $\tilde{K}$. But, as \ $x$ \ is a good point, for the given \ $\tilde{K}$ \  there exists a neighbourhood \ $V$ \ of \ $x$ \ and a compact set \ $L \subset G$ \ as in the definition. We may assume that \ $x_{\nu}$ \ is in \ $V$ \ for \ $\nu \geq \nu_0$ \ and so we may find, for \ $\nu \geq \nu_0, \alpha \geq \alpha(\nu)$, elements  \ $\gamma_{\nu, \alpha} \in L$ \ such that \ $\lim_{\alpha \to \infty} \gamma_{\nu,\alpha}.x_{\nu} = y_{\nu} \quad \forall \nu \geq \nu_0$.\\
Up to pass to a subsequence for each given $\nu \geq \nu_{0}$, we may assume that the sequence\ $(\gamma_{\nu, \alpha})$ \ converges to some \ $\gamma_{\nu} \in L$.  And again, that the sequence $\gamma_{\nu}$ converges to some $\gamma \in L$. So the continuity of \ $f$ \ gives  \ $y_{\nu} = \gamma_{\nu}.x_{\nu} \to \gamma.x = y \in \tilde{K}$ \ giving a contradiction because we assume \ $(G.x) \cap \tilde{K} = \emptyset$. This proves ii).\\

Let \ $E : = (U,B, j)$ \ be a \ $n-$scale on \ $\Omega$ \ adapted to the $n-$cycle \ $G.x$. Then the compact set \ $K : = j^{-1}(\bar U\times \partial B)$ \ does not meet \ $G.x$, by definition of an adapted scale. Using ii), there exists a neighbourhood \ $V$ \ of \ $x$ \ such that for any \ $x' \in V$ \ we have \ $(\overline{G.x' })\cap K = \emptyset$. As for a good point \ $x' \in V$ \ we know that \ $G.x'$ \ is a closed analytic subset, the \ $n-$scale is then adapted to \ $G.x'$. This implies that the dimension of \ $G.x'$ \ is at most equal to \ $n$. But the semi-continuity of the dimension of the stabilizers implies that the dimension of \ $G.x' \simeq G\big/St(x')$ \ is at least equal to \ $n = \dim (G\big/St(x))$. This proves iii).\\

Remark that for any $x' \in V$ such that $\overline{G.x'}$ is a closed analytic subset in $\Omega$, the previous  proof  shows also that $\overline{G.x'}$ is of pure dimension $n$.\\

To prove iv)  fix a good connected open set \ $\Omega$ \ and define
$$ Z : = \{ (g,x,y) \in G\times \Omega \times \Omega  \ / \  y = g.x .\}  $$
This is a closed analytic subset in $G\times \Omega\times \Omega$. Let us show that the projection \ $p : Z \to \Omega \times \Omega$ \ is semi-proper. Pick a point \ $(x,y) \in \Omega\times \Omega$ \ and choose compact neighbourhoods  \ $V$ \ and \ $K$ \ respectively of \ $x$  and of \ $y$ \ in \ $ \Omega$. Using the fact that any open set  \ $\Omega' \subset\subset \Omega$ \ is  uniformely good, for the compact \ $K$ \ in \ $\Omega$  we find a compact set \ $L$ \ in \ $G$ \ such that for \ $x' \in \Omega'$ \ and \ $g \in G$ \ with \ $g.x' \in K$ \ there exists \ $\gamma \in L$ \ with \ $\gamma.x' = g.x'$. Choosing \ $\Omega'$ \ containing \ $V$ this implies  that we have
$$ p(Z) \cap (V \times K) = p\big((L \times V \times K) \cap Z\big). $$
So the projection \ $p$ \ is semi-proper. Its image \ $p(Z)$ \ is then a closed analytic subset of \ $\Omega \times \Omega$ \ by Kuhlmann's theorem [K.64], [K.66]. But now the projection
$$ \pi : p(Z) \to \Omega $$
is \ $n-$equidimensional, thanks to iii), and has irreducible generic fibers on a normal basis \ $\Omega$. So its fibers (with generic multiplicity equal to $1$) define an analytic family of \ $n-$cycles of \ $X$ \ parametrized by \ $\Omega$. It is clearly \ $f-$ analytic because each fiber is irreducible{\footnote{But some multiplicities may occur.} and we have an holomorphic section because each \ $x $ \ lies in \ $G.x$.\\
To prove v) let us prove that the holomorphic map $\varphi : \Omega \to \mathcal{C}_{n}^{f}(\Omega)$ classifying the fibers of $p(Z)$ is semi-proper.  Fix a point $C \in \mathcal{C}_{n}^{f}(\Omega)$ and a point $x_{1}, \dots, x_{k} $ in each irreducible component of $\vert C\vert$. Let $W$ a relatively compact open neighbourhood of $\{x_{1}, \dots , x_{k}\}$ in $\Omega$ and let $\mathcal{W}$ be the open set in $\mathcal{C}_{n}^{f}(\Omega)$ of cycles such that each irreducible component meets $W$. Let  $C'$ be in $\mathcal{W} \cap \varphi(\Omega)$ we know that if $C' = \varphi(z)$ that  $\vert C'\vert = G.z$. So $G.z$ has to meet $W$ and we can choose $y $ in the compact set $ \bar W$ such that $C' =  \varphi(y)$  and this gives the semi-properness of $\varphi$. Now the semi-proper direct image theorem 2.3.2 of [B.15] implies that the image $Q$  of $\varphi$ is a finite dimensional complex space. \\
Then the holomorphic map $q : \Omega \to Q$ which is a quasi-proper equidimensional holomorphic quotient for the action $f$ on $\Omega$ as {\bf each} fiber of $q$ is set-theoretically  a $G-$orbit.$\hfill\blacksquare$\\

\subsection{The conditions ${\rm [H.1], [H.2]} $ and ${\rm [H.3]}.$}

Now we shall consider the following conditions on the action \ $f$.

\begin{itemize}
\item{\bf  There exists a $G-$invariant  dense open set $\Omega_{1}$  in $ X$  which admits a GF- holomorphic quotient}. \hfill ${\bf [H.1]}$
\end{itemize}

Recall that this means that there exists a $G-$invariant  geometrically f-flat holomorphic map $q : \Omega_{1} \to Q_{1}$ onto a reduced complex space $Q_{1}$  such that each fiber of $q$ over a point in $Q_{1}$ is set-theoretically an orbit in $\Omega_{1}$.\\

The following  stronger condition will be useful in the sequel:

\begin{itemize}
\item{\bf  There exists a $G-$invariant  dense  open set $\Omega$  in $ X$  which is   good for the action    $f$  (on $\Omega$)} .\hfill ${\bf [H.1str]}$
\end{itemize}

An immediate consequence of the proposition \ref{Good 2} is that the condition ${ \rm[H.1str]}$  below is sufficient to satisfy ${\rm [H.1}]$ by taking $\Omega_{1}$ as the set of normal points in $\Omega$ (which is open dense and $G-$invariant in $\Omega$). \\

Now assume [H.1] and  define  $\mathcal{R} : = \{ (x,y) \in \Omega_{1} \times \Omega_{1} \ / \  y \in G.x\} $. It is a closed analytic set in \ $\Omega_{1}\times \Omega_{1}$. We shall denote \ $\overline{\mathcal{R}}$ \ the closure of \ $\mathcal{R}$ \ in \ $X \times X$. Our second assumption will be :
\begin{itemize}
\item  {\bf The  subset $\overline{\mathcal{R}}$   is analytic  in  $X \times X$ and there exists an open dense subset $\Omega_{0} \subset \Omega_{1}$ such that for each $x \in \Omega_{0}$ we have} 
\begin{equation*}
{\bf \overline{G.x} = \overline{\mathcal{R}} \cap(\{x\}\times X)}. \qquad  \qquad \qquad \qquad \qquad \qquad  \qquad \qquad \qquad \qquad \qquad  {\bf [H.2]}
\end{equation*}
\end{itemize}

Remark that the first projection $p_{1} : \overline{\mathcal{R}}\cap (\Omega_{0}\times X) \to \Omega_{0}$ is quasi-proper because we have a holomorphic section of this map (with irreducible fibers, thanks to $[H.2]$) which is given by $x \mapsto (x, x)$.\\
Assuming that $\Omega_{0}$ contains only normal points\footnote{This not restrictive, as we may always assume that $X \setminus \Omega_{0}$ contains the non normal points in $X$. We shall always assume that $\Omega_{0}$ is normal in the sequel, without any more comment.} in $X$, the equidimensionality and quasi-properness  on $\Omega_{0}$ of the projection of $\overline{\mathcal{R}}$ implies that there exists a holomorphic map 
$$\bar \varphi_{0} : \Omega_{0} \to \mathcal{C}_{n}^{f}(X)  $$
where the supports are  given by $x \mapsto \overline{G.x}$ and where the multiplicity is generically equal to $1$.\\
We shall denote \ $\Gamma \subset \Omega_{0} \times \mathcal{C}_n^f(X) $ \ the graph of the map \ $\bar \varphi_{0}$ \ and \ $\bar \Gamma$ \ its closure in \ $X \times \mathcal{C}_n^f(X) $. We shall denote \ $\theta : \bar \Gamma \to X$ \ the map induced by the first projection. Our last hypothesis is :
\begin{itemize}
\item {\bf The  map  $\theta : \bar\Gamma \to X$  is  proper} . \hfill ${\bf [H.3]}$
\end{itemize}

This hypothesis ${\rm [H.3]}$, assuming ${\rm [H.1]}$ and ${\rm [H.2]}$ is in fact equivalent to ask that the projection $\overline{\mathcal{R}} \to X$ is strongly quasi-proper: this is an immediate consequence of the  proposition 3.2.2 of [B.15].\\

The following proposition shows that these conditions ${\rm [H.1], [H.2]}$ and ${\rm [H.3]}$ are necessary for the existence of a SQP-meromorphic quotient for a completely holomorphic action of $G$ on $X$.

\begin{prop}\label{necessary}
Assuming that the completely holomorphic action $f : G \times X \to X$ of the connected complex Lie group $G$ on the irreducible complex space $X$ has a SQP-meromorphic quotient, then the conditions ${\rm [H.1], [H.2]}$ and ${\rm [H.3]}$ are satisfied.
\end{prop}

\parag{Proof} The conditions to be a SQP-meromorphic quotient gives an open set $\Omega_{1}$ which is dense, $G-$stable and which admits a GF holomorphic quotient for the action on $\Omega_{1}$. So ${\rm [H.1]}$ is clear.\\
Let $S$ be the graph of the equivalence relation given by $q$ on $\tilde{X}$. Then the proper direct image $(\tau\times \tau)(S)$ is $\overline{\mathcal{R}}$ and so the condition ${\rm [H.2]}$ is satisfied as soon as we can find an open and dense subset $\Omega_{0}$ in $\Omega_{1}$ such that for each $x \in \Omega_{0}$ we have the equality 
$$\overline{G.x} = \overline{\mathcal{R}} \cap(\{x\}\times X). $$
But this property is given by the condition ii) in the definition of a SQP-meromorphic quotient.\\
The composition of $q$ with the classifying map $Q \to \mathcal{C}_{n}^{f}(\tilde{X})$ for the fibers of $q$ gives a holomorphic map $\tilde{\psi} : \tilde{X} \to \mathcal{C}_{n}^{f}(\tilde{X})$. Composed with the direct image map, which is holomorphic (see [B.M] ch.IV ; the ``quasi-proper'' part of this result is easy, as $\tau$ is proper) $\tau_{*} :\mathcal{C}_{n}^{f}(\tilde{X}) \to \mathcal{C}_{n}^{f}(X)$, we obtain the fact that $\overline{\mathcal{R}}$ is  strongly quasi-proper on $X$ which is the condition  ${\rm [H.3]}$.$\hfill \blacksquare$

 \subsection{Existence of SQP-meromorphic quotient}
 
Now we shall prove that  conditions \ ${\rm [H.1], [H.2]}$ \  and \ ${\rm [H.3]}$ \ on a completely holomorphic action of a connected complex Lie group $G$  on an irreducible complex space $X$  are sufficient for the existence of a SQP meromorphic quotient. 
 
 \begin{thm}\label{mero. quot.}
Under the hypothesis \ ${\rm [H.1], [H.2]}$ \  and \ ${\rm [H.3]}$ \ there exists a proper  $G-$equivariant  modification \ $\tau : \tilde{X} \to X$ \ with center contained in \ $X \setminus \Omega_{0}$\footnote{The dense open subset $\Omega_{0} \subset \Omega_{1}$ is defined in the condition ${\rm [H.2]}$} and a   geometrically f-flat holomorphic map
$$ q :  \tilde{X} \to  Q$$
on a  reduced complex space, which give a strongly quasi-proper meromorphic quotient for the given \ $G-$action.
\end{thm}

Of course the complex space \ $\tilde{X}$ \ is the topological space \ $\bar \Gamma$ \ with a structure of a reduced  complex space such that the projection on \ $X$ \ is a proper modification. Then the space \ $Q$ \ is the  image of \ $\tilde{X}$ \ in \ $ \mathcal{C}_{n}^{f}(X)$. So we need some semi-proper direct image theorem for such a map to prove this result. Such a result is the content of  the theorem 2.3.2 of [B.15]

\parag{Proof}  The first remark is that the hypothesis \ $[H. 3]$ \ says that the projection \ $p : \bar \Gamma \to X$ \ is a proper topological modification of \ $X$. But to apply directly the part ii)  of the theorem 2.3.6 of [B.13] to the projection \ $p_{1} : \overline{\mathcal{R}} \to X$ \ we need quasi-properness of this map. This is given by the proposition 3.2.2  of [B.15] as we have the condition $[H.3]$.\\
Then we obtain a proper (holomorphic) modification with center in \ $\Sigma \subset X \setminus \Omega_{0}$,  $\tau :  \tilde{X} \to X$,  and a \ $f-$analytic family of cycles in \ $X$ \   parametrized by \ $\tilde{X}$ \ extending the family \ $(\overline{G.x})_{x \in \Omega_{0}}$, corresponding to a ``holomorphic'' map extending \ $\bar \varphi_{0}$ :
$$\tilde{\varphi} : \tilde{X} \to \mathcal{C}_{n}^{f}(X) .$$
Now let us prove that this map \ $\tilde{\varphi}$ \ is quasi-proper\footnote{This makes sense as the fibers are closed analytic subsets of  \ $\tilde{X}$.}. This will allow us to apply the theorem 2.3.2  of {\it loc. cit.} \ and to define the reduced complex space \ $Q $ as the image $ \tilde{\varphi}(\tilde{X})$. Then it will be easy to check that the map \ $\tilde{\varphi} : \tilde{X} \to Q$ \ is a strongly quasi-proper  meromorphic quotient for the \ $G-$action we consider.\\
If \ $C_{0}$ \ is in \ $\mathcal{C}_{n}^{f}(X)$ \ and  is not the empty cycle, choose a relatively compact open set \ $W$ \ in \ $X$ \ such that any irreducible component of \ $\vert C_{0}\vert$ \ meets \ $W$. Then let \ $\mathcal{W}$ \ be the open set in \ $\mathcal{C}_{n}^{f}(X)$ \ defined by the condition that any irreducible component of \ $C$ \ meets \ $W$ \ for \ $C \in \mathcal{W}$. Then we shall prove that there exists a compact set \ $K$ \ in \ $\tilde{X}$ \ such that any irreducible component of the fiber of $\tilde{\varphi}$ at a point in $\mathcal{W} \cap \tilde{\varphi}(\tilde{X})$ meets $K$.  Let \ $K : = \tau^{-1}(\bar W)$. If \ $(y, C)$ \ is in \ $\tilde{X}$\footnote{Recall that, as a topological space, \ $\tilde{X} = \bar \Gamma$.} with \ $C \in  \mathcal{W}$, each irreducible component of \ $C$ \ meets \ $W$. But the fiber of \ $\tilde{\varphi}$ \ at \ $C$ \ is equal to \ $\vert C\vert$, and the quasi-properness is proved. $\hfill \blacksquare$\\

\section{Application.}

\subsection{The sub-analytic lemma.}

We shall use the following lemma (see  [G-M-O]) in our application.

\begin{lemma}\label{Bishop}
Let $M$ be a reduced complex space and $Y \subset M$ a closed analytic subset with no interior point in $M$. Let $R$   be a closed (complex)  analytic subset in  $ M \setminus Y$ such that $\bar R$ is a sub-analytic set in $M$. Then $\bar R$ is a (complex) analytic subset in $M$.
\end{lemma}

This important lemma is a consequence of Bishop's theorem (see [Bi.64]) and of a classical result on sub-analytic subsets (see [G-M-O]).

\subsection{Proof of the theorem \ref{G = K.B}.}

Now we shall assume that $G$ is a connected complex Lie group such that $G = K.B$ where $B$ is a closed  complex connected subgroup of $G$ and $K$ a compact real subgroup of $G$.

The first condition ${\rm [H.1str]}$ for the $G-$action is given by the following lemma :

\begin{lemma}\label{good B to good G}
In the situation of the  theorem \ref{G = K.B}, assume that we have a  \ $G-$invariant  open set $\Omega$ which is a good open set  for the $B-$action, then $\Omega$ is a good open set for the $G-$action.
\end{lemma}

\parag{Proof} Consider a point $x \in \Omega$ and a compact set $M$ in $\Omega$. Then there exists a neighbourhood $V$ of $x$ in $\Omega$ and a compact set $L$ in $B$ such that $b.y \in M$ for some $y \in V$ and some $b \in B$  implies that we can find $\beta \in L$ with $b.y = \beta.y$. Now assume that $M$ is $K-$invariant (here we use the $G-$invariance of $\Omega$) and that $g.y$ is in $M$ for some $g \in G$ \ and some $y \in V$. Write $g= k.b$ for some $k \in K$ and $b \in B$. Then $b.y$ is again in $M$ so we can find $\beta \in L$ with $\beta.y = b.y$ and then $ g.y = k.\beta.y$ with $k.\beta \in K.L$ which is a compact set in $G$.
So $x$ is a good point for the $G-$action on $\Omega$. $\hfill \blacksquare$\\

The corollary of the next lemma will give the first part of ${\rm [H.2]}$ for the $G-$action assuming that we have a  $G-$invariant  dense  good Zariski open set $\Omega$ for the $B-$action with the condition ${\rm [H.2]}$ for the $B-$action. In order to use the sub-analytic lemma \ref{Bishop} in this situation we shall need the following lemma.

\begin{lemma}\label{H.2 i)}
Let $\Omega$ be an open  $G-$invariant good set for the $B-$action, and then also good for the $G-$action thanks to the previous lemma.\\
 Let $\chi : K\times X\times X \to X\times X$ the map given by  $(k,x,y) \mapsto (k.x,y)$ and let \\
  $p : K\times X\times X \to X\times X$ be the natural  projection. Then we have
$$ p(\chi^{-1}(\overline{\mathcal{R}_{B}})) = \overline{\mathcal{R}_{G}} $$
where we define
 $$\mathcal{R}_{B} := \{(x,y) \in \Omega\times \Omega \ / \  B.x = B.y\} \quad {\rm and} \quad \mathcal{R}_{G} := \{(x,y) \in \Omega\times \Omega \ / \  G.x = G.y \} ,$$
 and where the closures are taken in $X\times X$.
 \end{lemma}
 
 \parag{Proof} Remark first that
  $$ p(\chi^{-1}(\mathcal{R}_{B})) = \{(x,y) \in \Omega\times \Omega \ / \  \exists k \in K \quadÊ B.k.x = B.y\}.$$
So $(x, y) \in p(\chi^{-1}(\mathcal{R}_{B}))$ implies   $y \in B.k.x \subset G.x$ and also $k.x \in B.y$;  we conclude that $x$ is in $K.B.y = G.y$.
 This gives the inclusion $ p(\chi^{-1}(\mathcal{R}_{B}))  \subset \mathcal{R}_{G}$. The opposite inclusion is easy because $G.x = G.y$ implies that $x \in K.B.y$ so there exists $k \in K$ such that $k.x \in B.y$.\\
 Now the map $\chi$ and  $p$ are proper, so we obtain the inclusion $ p(\chi^{-1}(\overline{\mathcal{R}_{B}})) \subset \overline{\mathcal{R}_{G}}$.\\
 Let $(x,y) := lim_{\nu\to \infty} (x_{\nu}, y_{\nu})$ where $(x_{\nu}, y_{\nu})$ is in $\mathcal{R}_{G}$ for each $\nu \in \mathbb{N}$. Then for each $\nu$  the exists $k_{\nu} \in K, b_{\nu} \in B$ such that $x_{\nu} = k_{\nu}.b_{\nu}.y_{\nu}$; then $(k_{\nu}^{-1}, x_{\nu}, y_{\nu}) \in \mathcal{R}_{B}$. Up to pass to a subsequence we may assume that  the sequence $(k_{\nu})$ converges to some $k \in K$. As $(k_{\nu}^{-1},x_{\nu},y_{\nu})$ is in $\chi^{-1}(\mathcal{R}_{B})$ for each $\nu$, the point $(k^{-1},x,y)$ is in 
  $$\overline{\chi^{-1}(\mathcal{R}_{B})} = \chi^{-1}(\overline{\mathcal{R}_{B}})$$
  and so $(x,y)$ is in $p(\chi^{-1}(\overline{\mathcal{R}_{B}})) $ proving the opposite inclusion.$\hfill \blacksquare$\\
 
 \begin{cor}\label{H.2 ii)}
 In the situation of the previous lemma, assume that $X \setminus \Omega$ is a (complex) analytic subset;  then  if the subset $\overline{\mathcal{R}_{B}}$ is (complex)  analytic in $X\times X$, the subset $\overline{\mathcal{R}_{G}}$ is also  a (complex) analytic subset of $X \times X$.
 \end{cor}
 
 \parag{Proof} Note first that the maps $\chi$ and $p$ are real analytic, so assuming that $\overline{\mathcal{R}_{B}}$ is analytic implies that $ p(\chi^{-1}(\overline{\mathcal{R}_{B}}))$ is sub-analytic. Then, as we know that $\mathcal{R}_{G}$ is an irreducible locally closed complex analytic subset, the conclusion follows from the lemma \ref{Bishop}, as our assumption that $\Omega$ is a Zariski (dense) open set in $X$ implies that $\Omega\times \Omega$ is Zariski open (and dense) in $X\times X$.$\hfill \blacksquare$\\
 
A first step to prove the quasi-properness of $\overline{\mathcal{R}_{G}}$ is our next result.
 
 \begin{lemma}\label{adherence orbites}
 Let assume that the $B-$action $f$ on $X$ satisfies ${\rm [H.1]}$ and ${\rm [H.2]}$. Let $\Omega_{0} \subset \Omega_{1}$ be an open set on which the fiber at any $x \in \Omega_{0}$  of $\overline{\mathcal{R}_{B}}$ is equal to $\overline{B.x}$ (with some multiplicity). Then the fiber at any $x \in \Omega_{0}$ of $\overline{\mathcal{R}_{G}}$ is equal to $\overline{G.x}$ (with some multiplicity).
 \end{lemma}
 
 \parag{Proof} As we know that the map  $x \mapsto \overline{B.x}$, with generic multiplicity $1$, extends to is a f-analytic family of cycles of $X$ parametrized by $\Omega_{0}$, for each sequence $(x_{\nu})_{\nu \in \mathbb{N}}$ of points in $\Omega_{0}$ converging to a point $x \in \Omega_{0}$ we have (with suitable multiplicity)  $\overline{B.x} = \lim_{\nu \to \infty} \overline{B.x_{\nu}}$ in the topology of $\mathcal{C}_{d}^{f}(X)$. We shall show that this implies, also with suitable multiplicity, the equality $\overline{G.x} = \lim_{\nu \to \infty} \overline{G.x_{\nu}}$ in the topology of $\mathcal{C}_{n}^{f}(X)$. As we have $G = K.B$ with $K$ compact, for any $y \in X$ we have $\overline{G.y} = K.\overline{B.y}$. So the inclusion of $\lim_{\nu \to \infty} \overline{G.x_{\nu}}$ in the fiber at $x$ of $\overline{\mathcal{R}_{G}}$ is clear. The point is to prove the opposite inclusion. Let $y$ be a point in the fiber at $x \in \Omega_{0}$ of $\overline{\mathcal{R}_{G}}$. It is a limit of a sequence $y_{\nu} \in G.x_{\nu}$ where $x_{\nu} \in \Omega_{0}$ converges to $x$. Write $y_{\nu} = k_{\nu}.b_{\nu}.x_{\nu}$ with $k_{\nu} \in K$ and $b_{\nu} \in B$. Up to pass to a subsequence, we may assume that the sequence  $(k_{\nu})$ converges to a point $k \in K$. So we have $k^{-1}.y$ which is the limit of the sequence $b_{\nu}.x_{\nu}$. We obtain that $k^{-1}.y$ is in the limit of $\overline{B.x_{\nu}}$ which has support equal to $\overline{B.x}$. Then $y$ is in $K.\overline{B.x} = \overline{G.x}$, concluding the proof.$\hfill \blacksquare$\\
 
 \parag{Proof of the theorem \ref{G = K.B}} The fact that the $G-$action satisfies ${\rm [H.1str]}$ is consequence of lemma \ref{good B to good G}. The analyticity of $\overline{\mathcal{R}_{G}}$ in $X \times X$ is proved at corollary \ref{H.2 ii)}. The lemma \ref{adherence orbites} gives a dense open set $\Omega_{0}$ where the fiber of the projection $p_{1}$ of  $\overline{\mathcal{R}_{G}}$ at each point $x \in \Omega_{0}$ is equal to $\overline{G.x}$ as a set.
 This implies the quasi-properness of $p_{1}$ over $\Omega_{0}$, because $x$ is in $G.x$ and $G$ is connected;  assuming (which is not restrictive) that $\Omega_{0}$ is normal, we obtain a holomorphic map 
 $$ \Phi : \Omega_{0} \longrightarrow \mathcal{C}_{n}^{f}(X)$$
 where the support of $\Phi(x)$ is equal to $\overline{G.x}$  for each $x \in \Omega_{0}$ and  with generic multiplicity equal to $1$. This complete the proof of ${\rm [H.2]}$ for the $G-$action.\\
 Thanks to proposition 3.2.2 of [B.15],  to prove ${\rm [H.3]}$  it is enough to show that the closure of the graph $\Gamma_{G}$ of $\Phi$ in $X\times \mathcal{C}_{n}^{f}(X)$ is proper on $X$.\\
  The projection $p_{B}: \overline{\mathcal{R}}_{B} \to X$ is strongly quasi-proper so, for $V$ compact in $X$, the limits of cycles (with convenient multiplicity) $\overline{B.x}$ for $x \in V \cap \Omega_{0}$ stay in a compact subset $S_{0}$ in $\mathcal{C}_{d}^{f}(X)$. Then let $S := K.S_{0}$; this is again a compact subset in $\mathcal{C}_{d}^{f}(X)$. But the subset $T$ of  limits of the generic fibers of the projection $p_{G}: \overline{\mathcal{R}}_{G} \to X$ for $x \in V$ is a compact set of $\mathcal{C}_{n}^{loc}(X)$ thanks to [B.13] theorem 2.3.6 i).
 For $x \in V \cap \Omega_{0}$ we have $p_{G}^{-1}(x) = \overline{G.x}$, and let $T'$ be the open dense set in $T$ described by the cycles $\bar \varphi_{0}(x), x \in V \cap \Omega_{0}$. As we know that for $x$ in $\Omega_{0}$ we have $\vert \bar \varphi_{0}(x) \vert = \overline{G.x}$ we know that each of them is an union of $d-$cycles in $S$, thanks to the lemma \ref{K-action} below. The proposition \ref{non inf. br.} implies that the set of  limits of $\overline{G.x_{\nu}}$ for $x_{\nu} \in V \cap \Omega_{0}$ is a compact set in $\mathcal{C}^{f}_{n}(X)$, proving [H.3].$\hfill \blacksquare$

\begin{lemma}\label{K-action}
Let $K$ be a compact Lie group with a holomorphic action on a reduced complex space $M$. Let $X$ be an irreducible complex space and let  $\varphi : X \to \mathcal{C}^{f}_{d}(M)$ and $\psi : X \to  \mathcal{C}^{f}_{n}(M)$ two holomorphic maps. Assume that for each $x$ in  a dense open subset $\Omega_{0}$ in $X$ we have
\begin{equation*}
 \vert \psi(x)\vert = \cup_{k\in K} \ k.\vert \varphi(x)\vert = K.\vert \varphi(x)\vert  \tag{@@@}
 \end{equation*}
where $K$ acts holomorphically on $\mathcal{C}^{f}_{d}(M)$ via the action of $K$ on $M$.\\
Then the relation $(@@@)$ holds for each $x \in X$.
\end{lemma}

\parag{proof} Let $x \in X$ and choose a sequence $(x_{\nu})_{\nu \geq 0}$ in $\Omega_{0}$ converging to $x$ and let $y$ a point in $\vert \psi(x)\vert$. Then we can choose a sequence $(y_{\nu})_{\nu \geq 0}$ of points in $\vert \psi(x_{\nu})\vert$ converging to $y$. As $x_{\nu}$ is in $\Omega_{0}$ we can write $ y_{\nu} = k_{\nu}.z_{\nu}$ with $z_{\nu}\in \vert \varphi(x_{\nu})\vert$ and $k_{\nu} \in K$. Up to pass to a subsequence we can assume that the sequence $(k_{\nu})$ converges to some $k \in K$. So the sequence $(z_{\nu})$ converges to $k^{-1}.y$ which is in $\vert \varphi(x)\vert$. This gives that $y$ is in $k.\vert \varphi(x)\vert$ and we have proved the inclusion $\vert \psi(x)\vert \subset K.\vert \varphi(x)\vert$.\\
Conversely, if $z$ is in $\vert \varphi(x)\vert$ and $k$ is in $K$, write again $z = \lim z_{\nu}$ with $z_{\nu} \in \vert \varphi(x_{\nu})\vert$ where the sequence $(x_{\nu})$ of points in $\Omega_{0}$ converges to $x$. We have $k.z_{\nu}\in \vert \psi(x_{\nu})\vert$, and then $k.z$ is in $\vert\psi(x)\vert$ conluding the proof.$\hfill \blacksquare$\\

\newpage

 \section{Bibliography}
 
 \begin{itemize}
 
  \item{[B.75]} Barlet, D. {\it Espace analytique r\'{e}duit des cycles analytiques complexes compacts d'un espace analytique complexe de dimension finie}, Fonctions de plusieurs variables complexes, II (S\'{e}minaire Fran\c cois Norguet, 1974-1975), pp.1-158. Lecture Notes in Math., Vol. 482, Springer, Berlin, 1975

 \item{[B.78]} Barlet, D. \textit{Majoration du volume des fibres g\'en\'eriques et forme g\'eom\'{e}trique du th\'eor\`{e}me d'aplatissement}, Seminaire Lelong-Skoda 78-79, Lect. Notes in Math. 822, p.1-17, Springer-Verlag, Berlin.
 
  \item{[B.08]} Barlet, D. \textit{Reparam\'etrisation universelle de familles f-analytiques de cycles et f-aplatissement g\'eom\'etrique} Comment. Math. Helv. 83 (2008),\\ p. 869-888.

 \item{[B.13]} Barlet, D. \ {\it Quasi-proper meromorphic equivalence relations}, Math. Z. (2013), vol. 273, p. 461-484.
 
 \itemÔ[B.15] Barlet, D. \ {\it Strongly quasi-proper maps and f-flattning theorem},\\ arXiv.1504.01579v1 [math. CV].
 
 \item{[Bi.64]} Bishop, E. {Condition for the analyticity of certain sets}, Michigan Math. J. 11 (1964), p.289-304.
 
 \item{[B-M]} Barlet, D. et Magn\'usson,J. {\it Cycles analytiques complexes I : Th\'eor\`{e}mes de pr\'eparation des cycles}, Cours Sp\'ecialis\'es 22, Soci\'et\'e math\'ematique de France, (2014).
 
  \item{[K.64]} Kuhlmann, N. \textit{{\"U}ber holomorphe Abbildungen komplexer R{\"a}ume} Archiv der Math. 15 (1964), p.81-90.
 
 \item{[K.66]} Kuhlmann, N. \textit{Bemerkungen {\"u}ber  holomorphe Abbildungen komplexer R{\"a}ume} Wiss. Abh. Arbeitsgemeinschaft Nordrhein-Westfalen 33, Festschr. Ged{\"a}achtnisfeier K. Weierstrass (1966), p.495-522.
 
 \item{[G-M-O]}  B. Gilligan, C. Miebach and K. Oeljeklaus {\it Homoheneous K{\"a}hler and Halmitonian manifolds}, Math. Ann. (2011), p.889-901.
 
  \item{[M.00]} Mathieu, D. \textit{Universal reparametrization of a family of cycles : a new approach to meromorphic equivalence relations}, Ann. Inst. Fourier (Grenoble) t. 50, fasc.4 (2000) p.1155-1189.
  
 \item{[R.57]} Remmert, R. \textit{Holomorphe und meromorphe Abbildungen komplexer R{\"a}ume}, Math. Ann. 133 (1957), p.328-370.

 \end{itemize}

\end{document}